\documentclass[a4paper,12pt]{article}
\usepackage[english]{babel}
\usepackage[T2A]{fontenc}
\usepackage[cp1251]{inputenc}
\usepackage[english]{babel}
\usepackage{amsthm}
\usepackage[tbtags]{amsmath}
\usepackage{amsfonts,amssymb}
\begin{document}

\newtheorem{theorem}{Theorem}
\newtheorem{lemma}{Lemma}
\newtheorem{proposition}{Proposition}
\newtheorem{Cor}{Corollary}

\begin{center}
{\large\bf Extensions of Centrally Essential Rings}
\end{center}

\begin{center}
Oleg Lyubimtsev\footnote{Nizhny Novgorod State University, Nizhny Novgorod, Russia; email: oleg\_lyubimcev@mail.ru .},
Askar Tuganbaev\footnote{National Research University MPEI, Moscow, Russia; Lomonosov Moscow State University, Moscow, Russia; tuganbaev@gmail.com .}
\end{center}
\textbf{Abstract.}
 A nonzero unital ring $R$ is said to be centrally essential if for every its nonzero element $a$, there exist nonzero central elements $x$ and $y$ with $ax = y$. In the paper, almost fully prime centrally essential rings are described in terms of ideal extensions, centrally essential Dorroh extensions, and trivial extensions. 


The work of Oleg Lyubimtsev is supported by Ministry of Education and Science of the Russian Federation, project FSWR-2023-0034, and scientific and educational mathematical center "Mathematics of technologies of the future". The study of Askar Tuganbaev is supported by Russian Science Foundation.

\textbf{MSC2020 database 16D25, 16R99}

\section{Introduction}\label{sec1}

All rings, considered in the paper, are associative nonzero unital rings, unless otherwise stated. A ring $R$ is said to be \textbf{centrally essential} if for any its nonzero element $a$, there exist nonzero central elements $x$ and $y$ such that $ax = y$. \footnote{It is clear that a ring $R$ with center $Z$ is  centrally essential if and only if the module $R_{Z}$ is an essential extension of the module $Z_{Z}$.} The class of centrally essential колец is a proper extension of the class of commutative колец. Centrally essential rings have been studied in \cite{MT19b}, \cite{MT20b}, \cite{LT24} and others. In  \cite[Proposition 2.8]{MT19b}, it is proved that semiprime centrally essential rings are commutative. 

For two rings $I$ and  $Q$, an \textbf{ideal extension} of $I$ by $Q$ is a  ring $R$ containing an ideal $I'$ such that $I'$  is ring-isomorphic to $I$ and $R/I'\cong Q$. For example, ideal extensions were studied in \cite{J16}. A ring is said to be \textbf{fully prime} (resp., \textbf{almost fully prime}) if every its ideal (resp., a nonzero proper ideal)  is prime; see \cite{Tsu94}, \cite{Tsu96}. We give an example of a noncommutative almost fully prime centrally essential ring $S$ which is an extension of its central ideal $I$ by the field $F$.

We recall that a \textbf{derivation} from  a ring $R$ into an $(R, R)$-bimodule $M$ is a mapping $d\colon R\to M$ such that 
$$
d(x) + d(y) = d(x + y) \text{ and } d(xy) = x\cdot d(y) + d(x)\cdot y \text{  for all }  x, y\in R.
$$ 
If $M$ is an $(R, R)$-bimodule with two derivations $d_1, d_2\colon R\to M$, then $d_1$ and $d_2$ are said to be \textbf{incomparable} if their kernels are not comparable with respect to inclusion, i.e. if there exist $x_1, x_2\in R$ such that $d_1(x_1)\neq 0$, 
$d_1(x_2) = 0$, $d_2(x_1) = 0$, $d_2(x_2)\neq 0$; see \cite{J16}.

\textbf{1.1. Example.} Let $R$ be a ring and let $F$ be a subring of $R$ which is a field. We consider two incomparable derivations $d_1, d_2\colon F\to R$ such that $d_1(a)d_2(b)\neq 0$ provided $d_1(a), d_2(b)\neq 0$. We define a ring $S$ of matrices 
$$
\left\lbrace\left.\begin{pmatrix}
f&d_1(f)&g\\
0&f&d_2(f)\\
0&0&f
\end{pmatrix}\;\;\right|\;\; f\in F,\, g\in R \right\rbrace.
$$
Then the ring $S$ is not commutative and is an extension of the central ideal $ I = 
\left\lbrace\left.\begin{pmatrix}
0&0&g\\
0&0&0\\
0&0&0
\end{pmatrix}\;\;\right|\;\; f\in F,\, g\in R \right\rbrace
$
by the field $F$; see \cite[Proposition 7 and Corollary 8]{J16}. In addition, $R$ is an almost fully prime centrally essential ring (see Теорему 1.2 below). 

Let $R$ and $S$ be two rings. We assume that $S$ is an $(R, R)$-bimodule such that the action $R$ on $S$ is coordinated with multiplication in $S$, i.e. $r(xy) = (rx)y$, $(xr)y = x(ry)$, $(xy)r = x(yr)$ for all $r\in R$ and $x, y\in S$. Then $D(R, S) = R\oplus S$ is a ring with multiplication
$$
(r_1, m_1)(r_2, m_2) = (r_1r_2, r_1m_2 + m_1r_2 + m_1m_2), \,\, \mbox{where} \,\, r_1, r_2\in R, m_1, m_2\in S
$$
which is called the \textbf{Dorroh extension} $R$ by the $S$. If $S = M$ is only an $(R, R)$-bimodule, then we set $m_1m_2 = 0$ (for all $m_1, m_2\in M$) and obtain the \textbf{trivial extension} $T(R, M)$ of $R$ by $M$ with multiplication
$$
(r_1, m_1)(r_2, m_2) = (r_1r_2, r_1m_2 + m_1r_2).
$$
Then $T(R, M)\cong \left\lbrace\left.\begin{pmatrix}
r&m\\
0&r\\
\end{pmatrix}\right|\;\; r\in R, m\in M \right\rbrace$ and $T(R, R)\cong R[x]/(x^2)$, where $(x^2)$ is the ideal  generated by $x^2$. We set 
$$
C_M(R) = \{r\in R \, \mid \, mr = rm \quad \forall \, m\in M\}
$$ 
and
$$
C_R(M) = \{m\in M \, \mid \, mr = rm \quad \forall \, r\in R\}.
$$ 
It is easy to see that $C_M(R)$ and $C_R(M)$ are a subring and a submodule of $R$ and $M$, respectively. We define subsets $T$ and $N$ in $R$ and $M$:
$$
T = \{r\in R \, \mid \, \exists \, c\in C_M(R)\cap Z(R): \, 0\neq rc\in C_M(R)\cap Z(R)\},
$$
$$
N = \{m\in M \, \mid \, \exists \, c\in C_M(R)\cap Z(R): \, 0\neq mc\in C_R(M)\}.
$$

Main results of the paper are Theorems 1.2 and 1.3.

\textbf{1.2. Theorem.} Let $R$ be an almost fully prime noncommutative ring. The following conditions are equivalent.
\begin{enumerate}
\item[{\bf 1)}]
 $R$ is a centrally essential ring.
\item[{\bf 2)}]
 $R$ is a ring with unique nonzero proper (left, right) ideal $P$ which is contained in the center of the ring, and the quotient ring $R/P$ is a field.
\item[{\bf 3)}]
$R$ is an extension of the central ideal $I$ by the field $F$.
\item[{\bf 4)}]
$R$ is an extension of the central ideal $I$ by the field $F$ which has two incomparable derivations from $F$ into some extension field $F$.
\end{enumerate}

\textbf{1.3. Theorem.} 
\begin{enumerate}
\item[{\bf 1.}]
 Let $R$ be a commutative ring, $S$ be an $(R, R)$-bimodule, $_RS = S_R$, and let $D(R, S)$ be the Dorroh extension $R$ by the $S$. The ring $S$ is a centrally essential ring if and only if $D(R, S)$ is a centrally essential ring.
\item[{\bf 2.}]
The trivial extension $T(R, M)$ is centrally essential if and only if $N = M$ and $T = R$.
\end{enumerate}

We denote by $Z(R)$ and $D(R)$ the center and set of all zero-divisors of the ring $R$, respectively. We recall some definitions used in the paper. A module is said to be \textbf{subdirectly indecomposable} if it contains the least nonzero submodule. The ring $R$ is said to be \textbf{right (left) subdirectly indecomposable} if the intersection of all its nonzero right (left) ideals is not equal to zero. A ring $R$ is said to be \textbf{subdirectly indecomposable} if the intersection of all its nonzero ideals is not equal to zero. A module is said to be \textbf{uniform} if any two nonzero submodules have nonzero intersection. If $R$ is a ring and every ideal of $R$ is idempotent, then $R$ is said to be \textbf{fully idempotent}. A ring is said to be \textbf{invariant} if all its one-sided ideals are ideals.

\section{Proof of Theorem 1.2}

$1)\Rightarrow 2)$. Let $R$ be a centrally essential ring. First, we assume that the set of ideals of $R$ is not linearly ordered. Then $R$ is a fully idempotent ring; see \cite[Theorem 2.1]{Tsu96}. Since the ring $R$ does not have nonzero nilpotent ideals, $R$ is semiprime. Then the ring $R$ is commutative which contradicts to the assumption of the theorem.

Let $R$ be a non-prime ring with linearly ordered set of ideals. Then $R$ has a minimal ideal $P$ such that every ideal of $R$, which is not equal to $P$, is idempotent; see \cite[Theorem 2.2]{Tsu96}. Therefore, the ring $R$ is subdirectly indecomposable. It follows from \cite[Proposition 2.1]{LT24} that $R$ is left and right subdirectly indecomposable. In particular, $R$ is a right uniform ring. Then set $D = D(R)$ of all (left or right) zero-divisors of $R$ is an ideal and quotient ring $R/D$ is commutative; see \cite[Lemma 2.2(2.2.3, 2.2.4)]{MT20b}. 

We verify that $D\subset Z(R)$. Let $0\neq a\in D$. Since $R$ is a centrally essential ring, $0\neq ac\in Z(R)$ for some $c\in Z(R)$. Let us assume that $c\in P$. Since $a\in D$, we have that $ay = 0$ for some $0\neq y\in R$. Then $a(yR) = 0$. Therefore,  $aP = 0$, since $P\subseteq yR$. But then $ac = 0$. This is a contradiction. Therefore, $c\notin P$. Since $cR$ is an  idempotent ideal, we have that $c = c^2b$ for some $b\in R$. For any $r\in R$, we have that
$$
cr(1 - cb) = cr - rc^2b = cr - rc =0.
$$
Therefore,  $cR(1 - cb) = 0$. Since the ideal $P$ is prime and $c\notin P$, we have $1 - cb = p\in P$. Therefore, $cb = 1 - p$ is an invertible element в $R$, since $p^2 = 0$. Consequently, $c$ is invertible and $a\in Z(R)$. Consequently, $D\subset Z(R)$. 

Now we verify that every nonzero proper ideal of $R$ coincides with $P$. Let $\overline{R} = R/P$. It follows from \cite[Theorem 1.3]{Tsu94} implies that $Z(\overline{R})$ is a field or $Z(\overline{R}) = 0$. The second case is impossible, since $R$ is a unital ring. Since $D\subset Z(R)$, we have that 
$\overline{D}$ is a central ideal в $\overline{R}$. Then every nonzero element of $\overline{D}$ is invertible in $\overline{R}$. Therefore,  
$\overline{D} = \overline{R}$ which is impossible. Therefore, $\overline{D} = \overline{0}$, $D = P$ and $R/D = R/P = \overline{R}$ is a field.
It follows from \cite[Theorem 1.7(3)]{LT24} implies that $R$ is an invariant ring. Consequently, $P$ is the unique nonzero proper (left, right) ideal в $R$. 

The implications $2)\Rightarrow 1)$ and $3)\Rightarrow 1)$ follow from \cite[Example 1.7]{MT19b}.

The implications $2)\Rightarrow 3)$ and $4)\Rightarrow 3)$ are obvious.

$3)\Rightarrow 4)$. By \cite[Theorem B]{J16},  if there does not exist incomparable derivations from $F$ into some field extension $F$, then there does not exist a  noncommutative  ring $R$ and a central ideal $I$ such that $R/I\cong F$.~$\square$

\section{Proof of Theorem 1.3}

\textbf{1.} Let $R$ be a commutative ring, $S$ be a centrally essential ring, $K = D(R, S)$, and let $(r, s)\notin Z(K)$. Then $s\notin Z(S)$. Since $S$ centrally essential ring, there exist elements $c, d\in Z(S)$ such that $0\neq sc = d$. If $rc + sc\neq 0$, then
$$
(0, 0)\neq (r, s)(0, c) = (0, rc + sc)\in Z(K).
$$ 
Let $rc + sc = 0$. Since $(r, s)\notin Z(K)$, we have that $(r, s)(0, 1) = (0, r\cdot 1 + s)\neq (0, 0)$. If $ r\cdot 1 + s\notin Z(S)$, then $0\neq (r\cdot 1 + s)c' = d'\in Z(S)$ for some $c', d'\in Z(S)$. Therefore, 
$$
(0, 0)\neq (r, s)\left[(0, 1)(0, c')\right] = (0, d')\in Z(K).
$$
Conversely, let $s\in S\backslash Z(S)$ Since $K = D(R, S)$ is a centrally essential ring, 
$$
(0, 0)\neq (0, s)(r, c) = (0, sr + sc)\in Z(K),
$$
where $(r, c)\in Z(K)$. Then $c, sr + sc$ are central elements of the ring $S$. It follows from $0\neq s(r\cdot 1 + c)\in Z(S)$ that $S$ is a centrally essential ring.

\textbf{2.} Let $R$ be a ring, $M$ be an $(R, R)$-bimodule, and let the ring $T(R, M)$ be centrally essential. Let $0\neq m\in M$. Since the ring $T(R, M)$ is centrally essential, for a matrix 
$\left.\begin{pmatrix}
0&m\\
0&0\\
\end{pmatrix}\right.$, there exists a central matrix $\left.\begin{pmatrix}
c&m'\\
0&c\\
\end{pmatrix}\right.$ such that 
$$
0\neq \left.\begin{pmatrix}
0&m\\
0&0\\
\end{pmatrix}\right.\cdot \left.\begin{pmatrix}
c&m'\\
0&c\
\end{pmatrix}\right. = \left.\begin{pmatrix}
0&mc\\
0&0\\
\end{pmatrix}\right.\in Z(T(R, M)).
$$
Then $c\in C_M(R)$. Indeed, if $cm''\neq m''c$ for some $m''\in M$, then
$$
\left.\begin{pmatrix}
0&m''\\
0&0\\
\end{pmatrix}\right.\cdot \left.\begin{pmatrix}
c&m'\\
0&c\\
\end{pmatrix}\right. \neq
\left.\begin{pmatrix}
c&m'\\
0&c\\
\end{pmatrix}\right. \cdot
\left.\begin{pmatrix}
0&m''\\
0&0\\
\end{pmatrix}\right.;
$$
this is a contradiction, since $\left.\begin{pmatrix}
c&m'\\
0&c\\
\end{pmatrix}\right.\in Z(T(R, M))$.
If $c\notin Z(R)$, then $cr\neq rc$ for some $r\in R$. Then
$$
\left.\begin{pmatrix}
c&m'\\
0&c\\
\end{pmatrix}\right.\cdot \left.\begin{pmatrix}
r&0\\
0&r\\
\end{pmatrix}\right. \neq
\left.\begin{pmatrix}
r&0\\
0&r\\
\end{pmatrix}\right. \cdot
\left.\begin{pmatrix}
c&m'\\
0&c\\
\end{pmatrix}\right.,
$$
and $\left.\begin{pmatrix}
c&m'\\
0&c\\
\end{pmatrix}\right.\notin Z(T(R, M))$. Consequently, $m\in N$ and $M = N$.

Now let $0\neq r\in R$. Then 
$$
0\neq \left.\begin{pmatrix}
r&0\\
0&r\\
\end{pmatrix}\right.\cdot \left.\begin{pmatrix}
c&m'\\
0&c\
\end{pmatrix}\right. = \left.\begin{pmatrix}
rc&rm'\\
0&rc\\
\end{pmatrix}\right.
$$
for some $\left.\begin{pmatrix}
c&m'\\
0&c\\
\end{pmatrix}\right., \left.\begin{pmatrix}
rc&rm'\\
0&rc\\
\end{pmatrix}\right.\in Z(T(R, M))$. Similar to the above, we verify that $c\in C_M(R)\cap Z(R)$. If $m(cr)\neq (rc)m$ for some $m\in M$, then
$$
\left.\begin{pmatrix}
rc&rm'\\
0&rc\\
\end{pmatrix}\right.\cdot \left.\begin{pmatrix}
0&m\\
0&0\\
\end{pmatrix}\right. \neq
\left.\begin{pmatrix}
0&m\\
0&0\\
\end{pmatrix}\right. \cdot
\left.\begin{pmatrix}
rc&rm'\\
0&rc\\
\end{pmatrix}\right.,
$$
whence we have $\left.\begin{pmatrix}
rc&rm'\\
0&rc\\
\end{pmatrix}\right.\notin Z(T(R, M))$. Consequently, $r\in T$ and $R = T$. 

Conversely, let $0\neq\left.\begin{pmatrix}
r&m\\
0&r\\
\end{pmatrix}\right.\in T(R, M)$. By assumption, there exists an element $c\in C_M(R)\cap Z(R)$ such that $0\neq rc = r'\in C_M(R)\cap Z(R)$. Then 
$$
\left.\begin{pmatrix}
r&m\\
0&r\\
\end{pmatrix}\right.\cdot \left.\begin{pmatrix}
c&0\\
0&c\\
\end{pmatrix}\right. = 
\left.\begin{pmatrix}
rc&mc\\
0&rc\\
\end{pmatrix}\right. = \left.\begin{pmatrix}
r'&mc\\
0&r'\\
\end{pmatrix}\right..
$$
If $mc\notin C_R(M)$, then, by assumption, $0\neq (mc)c'\in C_R(M)$ for some $c'\in C_M(R)\cap Z(R)$. Then $cc'\in C_M(R)\cap Z(R)$ and
$$
0\neq \left.\begin{pmatrix}
r&m\\
0&r\\
\end{pmatrix}\right.\cdot \left.\begin{pmatrix}
cc'&m'\\
0&cc'\
\end{pmatrix}\right. = \left.\begin{pmatrix}
rcc'&mcc'\\
0&rcc'\\
\end{pmatrix}\right.,
$$
where $\left.\begin{pmatrix}
cc'&0\\
0&cc'\\
\end{pmatrix}\right.$, $\left.\begin{pmatrix}
rcc'&mcc'\\
0&rcc'\\
\end{pmatrix}\right.\in Z(T(R, M))$. Consequently, $T(R, M)$ is a centrally essential ring.~$\square$

A not necessarily unital ring $R$ is said to be \textbf{centrally essential} if either it is commutative or for any its non-central element $a$, there exist nonzero central elements $x$ and $y$ such that $ax = y$. 

In Corollaries 3.1--3.3 below, the ring $R$ is not necessarily unital.

\textbf{3.1. Corollary.} A ring $R$ is centrally essential if and only if the ring $D(\mathbb{Z}, R)$  is centrally essential. Consequently, any centrally essential ring can be embedded in a centrally essential unital ring.

\textbf{3.2. Corollary.} If $T(R, M)$ is a centrally essential ring, then the ring $R$ is centrally essential.

\textbf{3.3. Corollary.} Let $I$ be an ideal of the ring $R$. The ring $T(R, I)$ is centrally essential if and only if the ring $R$ is centrally essential.

\section{Remarks and Additional Examples}

\textbf{4.1. Remark.} If the characteristic of the field $F$ from Theorem 1.2(3, 4) is equal to zero, then we can add an equivalent condition to formulation of the theorem, see \cite[Theorem 2]{J16}: 

\noindent {\bf 5)} The transcendence degree of the field $F$ over its prime subfield exceeds 1.

\textbf{4.2. Proposition.} Let $R$ be an almost fully prime centrally essential ring. Then the additive group $R^+$ of the ring $R$ is isomorphic to $\bigoplus_{\alpha} Q$ or $\bigoplus_{\alpha} Z_p$, for some prime integer $p$, where $\alpha$ is an arbitrary cardinal number.

\textbf{Proof.} Let $P$ be the unique	nonzero proper ideal of the ring $R$. It follows from Theorem 1.2 that $R/P$ is a field. Consequently,  $char (R/P) = 0$ or $char (R/P) = p$, where $p$ is a prime integer. Then $(R/P)^+$ is a vector space over the field $\mathbb{Q}$ or $\mathbb{Z}_{p}$. It follows from $P^2 = 0$  that $P$ is an $R/P$-module. Therefore, $P$ also is a vector space over $\mathbb{Q}$ or $\mathbb{Z}_{p}$. Since  $R^+\cong P^+\bigoplus (R/P)^+$, we have that $R^+\cong\bigoplus_{\alpha} Q$ or $R^+\cong\bigoplus_{\alpha} Z_p$.~$\square$

\textbf{4.3. Remark.} There exist almost fully prime rings having distinct nonzero proper ideals; see \cite[Example 2.4]{Tsu96}. Moreover, there exist fully prime noncommutative  rings having infinitely many ideals; see \cite{Tsu94}.

\textbf{4.4. Remark.} Let $D = D(R)$ be the set of all zero-divisors of the ring $R$. If $\{0\}\neq D\subset Z$, then $R$ is centrally essential ring.  It follows from \cite[Example 1.1]{KB98} that there exists a noncommutative  centrally essential ring such that all its zero-divisors are contained in the center of the ring. 

\textbf{4.5. Example.} Let $\mathbb{Z}$ be the ring of integers, $\mathbb{Q}$ be the field of rational numbers, and let $\mathbb{Q}[i]$ be the field of fractions of the ring of Gaussian integers $\mathbb{Z}[i]$. We denote by $\varphi$ the automorphism $q_1 + q_2i \mapsto q_1 - q_2i$ of the field $\mathbb{Q}[i]$. Let $R_X$ and $R_Y$ be the localizations $\mathbb{Z}[i]$ with respect to prime ideals $X = (2 + i)$ and $Y = (2 - i)$, respectively. We set $R = R_X\cap R_Y$. The ring $R$ coincides with the set of all irreducible rational Gaussian fractions whose denominators are divided neither by $2 + i$ nor by $2 - i$. We denote by $M$ the right $R$-module $\mathbb{Q}[i]/R_X$. We also turn $M$ into a left $R$-module by the rule $rm = m\varphi(r)$ for all elements $r\in R$ and $m\in M$. It is directly verified that $M$ is an $(R, R)$-bimodule. We consider the trivial extension $S = T(R, M)$ of the ring $R$ by the module $M$. Then $C_M(R) = \mathbb{Q}$. If $0\neq z = q_1 + q_2i$, where $q_2\neq 0$, then $z\notin T$. By Theorem 1.3(2), the ring $S$ is not centrally essential.

We give an example of a centrally essential trivial extension $S = T(R, M)$, where $_RM\neq M_R$.

\textbf{4.6. Example.} We construct the required extension with the use of Example 1.1 from \cite{KB98}. Let $\mathbb{Z}[x, y]$ be the polynomial ring in two variables $x$ and $y$ over $\mathbb{Z}$, $K = \mathbb{Z}[x, y]/(y^2)$, and let $R = K[t, \sigma]$ be the differential polynomial ring $\sum t^ia_i$, $a_i\in K$ with derivation $\sigma$ defined as $\sigma x = y$ and $\sigma y = 0$. Since $y^2 = 0$, we have that
$$
(f(x) + g(x)y)t = t(f(x) + g(x)y) + f'(x)y.
$$
We consider the ring 
$$
S = \left\lbrace\left.\begin{pmatrix}
f(x, y)& g(t)\\
0&f(x, y)\\
\end{pmatrix}\right|\;\; f(x, y)\in K, g(t)\in R \right\rbrace.
$$
Since
$$
\left.\begin{pmatrix}
x&0\\
0&x\\
\end{pmatrix}\right.\cdot \left.\begin{pmatrix}
0&t\\
0&0\\
\end{pmatrix}\right. = \left.\begin{pmatrix}
0&xt\\
0&0\\
\end{pmatrix}\right.\neq
\left.\begin{pmatrix}
0&tx\\
0&0\\
\end{pmatrix}\right. = 
\left.\begin{pmatrix}
0&t\\
0&0\\
\end{pmatrix}\right.\cdot \left.\begin{pmatrix}
x&0\\
0&x\\
\end{pmatrix}\right.,
$$
the ring $S$ is not commutative. In addition, $S$ is a centrally essential ring. Indeed, let $0\neq \left.\begin{pmatrix}
k& f(t)\\
0&k\\
\end{pmatrix}\right.\notin Z(S)$, $k = f(x) + g(x)y$. If $f(x) = 0$, then $g(x)y\in C_M(R)$. If $f(x)\neq 0$, then  $0\neq ky\in C_M(R)$, where 
$y\in C_M(R)$. Therefore, $T = R$. Further, if $f(t)\notin C_R(M)$, then $f(t)$ contains a monomial with nonzero coefficient $h(x)$ and 
$0\neq f(t)y\in C_R(M)$. Consequently, $N = M$ and $S$ is a centrally essential ring by Theorem 1.3(2).

\end{document}